\documentclass[a4paper, 1pt]{article}
\usepackage{amssymb,amsbsy,amsmath,amsfonts,amssymb,amscd}
\usepackage{color,latexsym,eucal,exscale,epic,eepic,epsfig}
\date{}
\numberwithin{equation}{subsection}
\newtheorem{theorem}{Theorem}[subsection]
\newtheorem{lemma}{Lemma}[subsection]
\newtheorem{proposition}[theorem]{Proposition}
\newtheorem{corollary}[theorem]{Corollary}
\newtheorem{definition}{Definition\rm}[subsection]

\def\ti{\widetilde }
\def\C{{\mathbb C}}
\def\NN{{\mathbb N}}

\def\R{{\mathbb R}}

\def\CC{{\mathcal C}}
\def\n{{\mathfrak N}}

\def\l{{\ell }}

\def\O{{\mathcal O}}

\def\L{{\mathcal L}}
\def\U{{\mathcal U}}
\def\W{{\mathcal W}}

\begin{document}

\title{HOLOMORPHIC EXTENSION OF DECOMPOSABLE DISTRIBUTIONS
 FROM A 
CR SUBMANIFOLD OF $\C^L$}

\author {Nicolas Eisen}

\maketitle

\abstract{ \noindent   Given $N$ a non generic smooth CR submanifold
of
$\C^L$, $N=\{(\n,h(\n))\}$ where $\n$ is generic in $\C^{L-n}$ and $h$ is a
CR map from $\n$ into $\C^n$. We prove, using only elementary tools,
 that if $h$ is decomposable at $p'\in \n$ then any decomposable 
CR distribution on
$N$ at $p=(p',h(p'))$
extends holomorphically to a complex transversal
wedge. This gives an elementary proof of the well known equivalent for
totally real non generic submanifolds, i.e if $N$ is a smooth totally
real
submanifold of $\C^L$ any continuous function on $N$ admits a
holomorphic
extension to a complex transverse wedge.}
\setcounter{section}{0}

\section{Introduction} 

\subsection {Statement of Results}
Let $\n$ be generic submanifold of $\C^{k+m}$ of CR dimension $k$
 and $h$ a  CR
map from $\n$ into some  $\C^{n}$ verifying $dh(0)=0$. Set $L=k+m+n$,
we construct a CR
submanifold $N$
of $\C^L$ near the origin as the graph of $h$ over $\n$, that is
$N=\{(\n,h(\n))\}$. It turns out that any non generic 
CR submanifold of $\C^L$
can be obtained in that fashion, see for example \cite{[Bo]}. The main 
question we address in this paper is the possible holomorphic
extension of a CR distribution of $N$ to some wedge $\W$ in a
 complex transverse direction.
 The aim of
 this paper is to give a proof of an extension result using only
 elementary tools. The CR structure of $N$ is determined by $\n$
hence any CR distribution on $N$ is a CR distribution on $\n$. 
\begin{definition}\label{D1} {\bf( a)} Let $\n$ be a smooth generic submanifold
of
$\C^L$
 A CR distribution $u$ on $\n$ is {\bf decomposable}
at
the point $p \in \n$ if, near $p$, $u=\sum_{j=1}^KU_j$ where the $U_j$
are 
CR distributions extending holomorphically
 to wedges $\W_j$ in $\C^{k+m}$  with edges
$\n$. We shall say a distribution $u$ on $N$ is decomposable at a
point $p=(p',h(p'))$ if $u$ is
decomposable at $p'$  on $\n$.

\medskip

 {\bf( b)} We say that $v$ is {\bf complex transversal} to $N$ at $p$
if
$v \not \in \C span T_pN$.
\end{definition}

Our main result is the following

\begin{theorem} \label{t1} Let $N=\{(\n,h(\n))\} $ 
be a non generic smooth ($\CC^{\infty}$) CR submanifold of $\C^{k+m+n}$ 
such that the map $h$ is decomposable at some $p_0'\in \n$.
Let  $v$ be a complex transversal vector to $N$ at
$p_0=(p_0',h(p_0'))$.
 If $f$ is a decomposable
CR distribution at $p_0\in N$ then, near $p_0$, there exists a wedge
$\W$ of direction $v$ whose edge contains a neighborhood of 
$p_0$ in $N$
and $F\in \O(\W)$ such that the boundary value of $F$ is $f$.
 Further more, there exist $\{F_{\l}\}_{\l=1}^n$, $F_{\l}\in
\O(\W)$ 
such that $dF_1\wedge ...\wedge dF_n\not = 0$ on $\W$ and each
$F_{\l}$ vanishes to order one on $N$.
\end{theorem}

The boundary value of a holomorphic function $F$ is defined by
$
\lim_{\lambda \to 0^+} \int_{\n}F(x+\lambda v)\varphi(x)dx,$ where $v
\in \W$. It turns
out
that if $F$ has slow growth in a wedge $\W$ (there exists a constant $C>0$ and a positive integer
$\ell$
such that
$$|F(z)|\leq {{C}\over{|dist(z,M)|^{\ell}}},$$
where $dist(z,M)$ denotes the distance from a point $z$ to $M$) then
the boundary value of $F$ defines a CR distribution on $N$ (see, for
example
\cite{[BER]}). We call
the integer $\ell$ above the {\bf growth degree } of $F$.

\bigskip

{\bf Remarks on the smoothness of $N$}

\bigskip

Note that one does not need $N$ to be smooth to be able to define a
 decomposable CR distribution on $N$. Indeed, suppose $F$ is
a holomorphic function of slow growth of growth degree $\ell$,  then one can
prove, following theorem 7.2.6 in \cite{[BER]} the next result.

\begin{proposition} Let $F$ be as above and suppose that the edge of
the
wedge $N$ is of regularity $\ell+1$, then the boundary value of $F$
defines
a CR distribution of order $\ell+1$ on $N$.
\end{proposition}

We thus define for $u$ a decomposable distribution, $u=\sum bv F_j$, 
the growth degree of $u$ to be the maximum of th growth degrees of the
$F_j$. We see that if the growth degree of a decomposable distribution $u$
is $\ell$ it makes sense to speak of a decomposable
distribution 
on a manifold of smoothness $\ell+1$.
 Hence the hypothesis of
smoothness
on $N$ in theorem\ref {t1}  can be replaced by $\ell +1$.

\bigskip

If instead of CR distribution we wish to consider functions, we can
reduce the smoothness hypothesis in theorem\ref {t1}.  A $\CC^0$
function $u$ 
that is decomposable
near a point $p$ is NOT the sum of $\CC^0$ functions $U_j$ extending
holomorphically. On the other hand if $u \in\CC^{\alpha},~\alpha \not
\in \NN$, if $u$ extends holomorphically at some $p$ into a wedge 
$\W_j$ of direction $w_j$ then the wedge $\W_j$ can be written as
$(\n \cap V_p)+i\Gamma_j$ where $V_p$ is a neighborhood of $p$ and
 $\Gamma_j$ is a conical neighborhood
of $w_j$ in the normal space to $N$ at $p$.
 The wave front set of $u$ (with respect to any micro-local
class) at $p$ is contained in the dual cone of $\Gamma_j$ denoted by
$\Gamma_j^0$. Hence if $u\in\CC^{\alpha}$ is decomposable at some $p$,
since the $\Gamma_j^0$ are pairwise disjoint the regularity
of the sum of the $U_j$ cannot be any better than the regularity of
each $U_j$. Therefore if $u$ is a  $\CC^{\alpha}$ decomposable
function, $u=\sum U_j$, then each $U_j$ has at least the same
regularity
as $U$. Hence if we wish to study the problem of holomorphic extension from the
point of view of functions rather than distributions we can replace
the smoothness of the manifold $N$ by $\CC^{1+\alpha}$ and study the
extension of   $\CC^{0+\alpha}$ 
decomposable CR
functions. We get the following result
\begin{theorem} \label{C} 
Let $N=\{(\n,h(\n))\} $ 
be a non generic $\CC^{1+\alpha}$ CR submanifold of $\C^{k+m+n}$ 
such that the map $h$ is decomposable at some $p_0'\in \n$.
Let  $v$ be a complex transversal vector to $N$ at
$p_0=(p_0',h(p_0'))$. If $f$ is a  $\CC^{0+\alpha}$ 
decomposable CR
function at $p_0$, then, near $p_0$, there exists a wedge
$\W$ of direction $v$ whose edge contains a neighborhood of 
$p_0$ in $N$
and $F\in \O(\W)$ such that
$F|_N=f$. Further more, there exist $\{F_{\l}\}_{\l=1}^n$, $F_{\l}\in
\O(\W)$ 
such that $dF_1\wedge ...\wedge dF_n\not = 0$ on $\W$ and each
$F_{\l}$ vanishes to order one on $N$.
\end{theorem}

We also obtain, using theorem \ref{t1} the following corollary

\begin{corollary}\label{D} Let $M$ be a  $\CC^{\infty}$ 
generic submanifold of $\C^L$
containing through some $p_0 \in M$ a proper  $\CC^{\infty}$   CR submanifold $N=(\n,h(\n))$ of
same CR dimension, $p_0=(p'_0,h(p'_0))$ with $p'_0\in \n$. Assume that
the function $h$ decomposable at $p_0$.
 Let 
$v\in T_{p_0}M \setminus 
\left [span_{\C}T_{p_0}N\right ]$. 
If $f$
is a CR distribution on $N$ that is  decomposable 
at $p_0$, then there exists a wedge
$\W$ in $M$
of direction $v$  whose edge contains a neighborhood of $p_0$ in $N$
and $F$ a $\CC^{\infty}$ CR function on $\W$ such that
$F|_N=f$. 
 Furthermore, there
exists a collection of $\CC^{\infty}$ CR functions $\{g_{\l}\}_{\l=1}^n$
vanishing to  order one  on $N$ and such that $dg_1\wedge....\wedge
dg_n \not=0 $ on $\W$.
\end{corollary}

The above corollary does not hold in the abstract CR structure case,
we finish the paper by constructing an abstract CR structure on which
there is no CR extension. More precisely, set $L={{\partial} \over {\partial \bar z}} +C(z,s,t){{\partial} \over {\partial s   }}
+tD(z,s,t){{\partial} \over {\partial t   }}$ and define
$L^0=L|_{t=0}$. The question we now address is the following: if
$f=f(z,s)$ is such that $L^0(f)=0$, does there exist $g$ such that
$L(f+tg)=0$? The answer in general is negative.

\bigskip

\begin{proposition} \label{P} There exist  $L$ as above and $h$ a real
analytic function, with $L^0(h)=0$ and such that the equation 
$L(h+tg)=0$ has no solution for $g \in {\mathcal C}^1$.
\end{proposition}

\subsection{Remarks} 
 
As pointed out in \cite{[Ei]}, this type of extension result was well
known
in the totally real case and is essentially due to Nagel \cite{[Na]}.
It can be stated in the following way:

\begin{theorem}\label{wk1} Let $N$ be a non generic totally real smooth
submanifold of $\C^L$. 
Let $ v \in \C^L$ be complex transversal to $N$ at $p$. Then for any
continuous function $f$ on a neighborhood of $p$ in $N$
 there exists $\W_v$, a wedge of direction
$v$  whose edge contains $N$
 such that $f$ has a holomorphic extension to $\W_v$.
\end{theorem}

This paper provides the easiest and simplest proof 
of this result since any continuous function on a totally real
submanifold is decomposable.

\bigskip 

We wish to point out the main differences between this paper and
\cite{[Ei]}.
In \cite{[Ei]} we obtain similar results, but the technics we use are
not at all the same and they yield extension results for non
decomposable CR distributions, also, the size of the wedges obtained are much  larger
than the one obtained here, roughly speaking, in \cite{[Ei]}  we
obtain
wedges that contain $(\n \times \R^n) \cap \{t_1>0\}$. 

\bigskip

As it is noted in \cite{[Tr-3]} on most CR submanifolds of $\C^L$
all CR functions are decomposable, hence the hypotheses of theorem
\ref{t1} hold in a generical sense for CR
distributions.
 However, they are examples of
CR submanifolds of $\C^L$ on which indecomposable CR functions exist.

\bigskip

Theorem \ref{t1} implies that the extension obtained 
is not unique, which differs greatly with the holomorphic extension
results obtained for generic submanifolds where the extension if it
exists is always unique.

\bigskip

One should note that the question of CR extension can be viewed as a
Cauchy problem with Cauchy data on a characteristic set $N$.

\bigskip

\subsection{Background}

For a general background on CR geometry, we refer the reader to the
books of Baouendi, Ebenfelt and Rothschild \cite{[BER]}, Boggess 
\cite{[Bo]} and Jacobowitz \cite{[Ja]}.

Most of the results on holomorphic extension deal with generic
submanifolds
of $\C^n$. In a general way, these results imply a forced unique
extension of CR functions under certain hypothesis on the manifold $M$
 such as Lewy non
degenerateness or more generally, minimality. Under these hypothesis
one does show that it is possible to fill a wedge with edge $M$ with
analytic discs attached to $M$. Using the fact that continuous CR
functions are uniform limits of polynomials and the maximum principle
one obtains a unique extension for continuous CR functions, see for
example the survey paper by Tr\'epreau \cite{[Tr-2]}. 

The subject of decomposable CR functions has been 
studied by many authors and it was believed that all CR were
decomposable, however Tr\'epreau produced examples of non decomposable
CR functions \cite{[Tr-3]} (an elementary explanation of this can be
found in a paper by Rosay \cite{[Ro]}). However one should note that
on most CR submanifolds of $\C^n$ any CR
function
is decomposable.

The subject of CR extension has not been studied in as much depth as
the
holomorphic extension has. When studying CR extension from a
submanifold of lower CR dimension, the tools involving analytic
discs still work, see for example \cite{[Tu-2]} 
and \cite{[Hi-Ta]}, however this is not
the case when the CR dimensions are equal see \cite{[Ei]}.

\section{Proof of the Extension Theorem}

\begin{proposition}\label{p1} Let $N=\{(\n,0)\} $ 
be a smooth  non generic CR submanifold of
$\C^{k+m+n}$. Let $\ti \W$ be a wedge in $\C^{k+m}$ with
edge $\n$ near $p'_0\in \n$. Suppose $F$ is a holomorphic function in
$\ti \W$ which is of slow growth, denote by $f$ its boundary value on
$\n$. For any $v$ complex transversal to $N$ at $p_0$, there exists a
wedge
$\W_v$ of direction $v$ whose edge contains $\n$ such that $f$ extends
holomorphically to $\W_v$.
\end{proposition}
\noindent{\bf Proof of Proposition \ref{p1}}. We begin with a choice of local coordinates on $\n$. 
$\n$ is a generic manifold in $\C^{k+m}$. We introduce local
coordinates near $p_0$. We may choose a local embedding so that
$p_0=0$ and 
$\n$ is parameterized in $\C^{k+m}=\C^k_z \times \C^m_{w'}$ by
 
\begin{equation} \label{e1} 
\n=\{(z,w')\in \C^k \times \C^m:Im(w')=a(z,
Re(w')),~~ a(0)=da(0)=0\}.
\end{equation}
We will denote by $s=Re(w') \in \R^m$, we thus have

\begin{equation} \label{e2} 
\n=\{(z,s+ia(z,s)\}\subset \C^k \times \C^m,~~~ T_0\n=\C^k\times \R^m.
\end{equation}

Define 
$\C T_p\n=T_p\n \otimes \C$ and
$T^{0,1}_p\n=T^{0,1}_p\C^{k+m} \cap \C T_p\n$. We say that $\n$ is a CR
manifold if dim$_{\C}T^{0,1}_p\n$ does not depend on $p$. The CR
vector fields of $\n$ are vector fields $L$ on $\n$ such that for any 
$p\in \n$ we have $L_p \in T^{0,1}_p\n$. One can
choose a basis $\L$ of $ T^{0,1}\n$ near the origin consisting of vector fields $L_j$
of the form

\begin{equation}\label{e3} 
L_j={{\partial} \over {\partial \overline z_j}}+\sum_{\ell=1}^nF_{j
\ell} {{\partial} \over {\partial  s_{\ell}}}.
\end{equation}

The wedge $\ti \W$ in $\C^{k+m}$ with edge $\n$ on which $F$ is
  defined
is given, in a neighborhood of the origin by 
$$
\ti \W=\left (\U+i\Gamma\right ),
$$
where $\U$ is a neighborhood of the origin in $\n$ and $\Gamma$ is
 a conic neighborhood of some vector $\mu$  in $\R^m\setminus
 \{0\}$. 
Note then that $F$ admits (trivially) a holomorphic extension
to the region $\ti \W\times \C^n\subset \C^{k+m+n}$.
This region is much more than a wedge. But it clearly contains
a wedge in $\C^{k+m+n}$ with direction $u$ whenever $u$ is a vector of the
$u=(u',u'')\in \C^{k+m}\times \C^n$ with $u'\in \ti \W$.
By \ref{e1} complex transversility of $v=(v',v'')$ means that
$v''\neq 0$.
Fix a vector $u\in \ti \W$. Consider a $\C$ linear change
of variables $T$ that is the identity on $\C^{k+m}\times
\{ 0\}$ and such that $T(v)=(u,v'')$. 
The desired extension of $f$
to a wedge of direction $v$ is then given by $F(T(z,w))$.
 We now need to show that the boundary value of $F(T(z,w))$ on $\n$ is
 $f$. The boundary value of $F$ on the wedge $\W$ is defined to
 be
\begin{equation}\label{ec1}
<f,\varphi>=\lim_{\lambda \to 0^+} \int_{\n} F(x+\lambda \gamma)\varphi(x)dx.
\end{equation}
For $\varphi \in \CC^{\infty}_0(\n)$ and $x+\lambda \gamma \in \ti
\W$.
Write $T=(T',T'') \in \C^{k+m}\times\C^{n}$. Let
$\tau=\tau(x,\lambda,\eta)$
be defined by 
$T'(x+\lambda \eta)=(x+\lambda\tau(x,\lambda,\eta))$. Since $T$ 
is the identity on $\C^{k+m}\times
\{ 0\}$ , we have $\lim_{\lambda\to 0^+}\lambda\tau(x,\lambda,\eta) =0$.
The boundary value of  $F(T(z,w))$ on $\n$ on the wedge $\W_v$ is then
given by
\begin{equation}\label{ec2}
\lim_{\lambda \to 0^+} \int_{\n} F(T(x+\lambda \eta))\varphi(x)dx,
\end{equation}
where $(x+\lambda \eta)\in \W_v$.
We then define 
$$G_{\tau}(\lambda)= \int_{\n}
F(x+\lambda\tau(x,\lambda,\eta))\varphi(x)dx,$$
$$F_{\gamma}(\lambda)=\int_{\n} F(x+\lambda \gamma)\varphi(x)dx.$$
Then by proposition 7.2.22 p189 of \cite{[BER]} we have
$$G_{\tau}(\lambda)-F_{\gamma}(\lambda)=O(\lambda), ~~~\lambda \to
0^+.$$
Hence the boundary value defined by \ref{ec2}
is equal to the boundary value defined by \ref{ec1}.
$\blacksquare$

\bigskip

We immediately note that if the boundary value is at least continuous
then the proof of the proposition yields the following
 
\begin{proposition}\label{p12}  Let $N=\{(\n,0)\} $ 
be a $\CC^1$  non generic CR submanifold of
$\C^{k+m+n}$. Let $\ti \W$ be a wedge in $\C^{k+m}$ with
edge $\n$ near $p'_0\in \n$. Suppose $F$ is a holomorphic function in
$\ti \W$ which has a continuous boundary value on $N$, denote by $f$ its boundary value on
$\n$. For any $v$ complex transversal to $N$ at $p_0$, there exists a
wedge
$\W_v$ of direction $v$ whose edge contains $\n$ such that $f$ extends
holomorphically to $\W_v$.
\end{proposition}

Using proposition \ref{p1}, we obtain a special case of theorem
\ref{t1}
from which we will deduce the later, namely.
\begin{proposition}\label{t'1} Let $N=\{(\n,0)\}$ be a  $\CC^{\infty}$ (resp
$\CC^{1+\alpha})$  CR submanifold of $\C^{k+m+n}$ .
Let  $v$ be a complex transversal vector to $N$ at $p_0$.
If $f$ is a decomposable distribution at $p_0$ 
(resp a $\CC^{0+\alpha}$  decomposable function ), then, near $(p_0,0)\in N$, there exists a wedge
$\W$ of direction $v$ whose edge contains a neighborhood of 
$(p_0,0)$ in $N$
and $F\in \O(\W)$ such that
$bv F=f$. 
\end{proposition}

\noindent {\bf Proof of Proposition \ref{t'1}.} 
Let $v$ be a complex transversal vector and let $u$ be a  CR
distribution
 on $\n$. By hypothesis,
$u=\sum_{j=1}^{K}
U_j$ where each $U_j$ is a boundary value of $F_j$, $F_j\in O(\ti
\W_j)$, $F_j$ of slow growth (or $F_j \in  \CC^{0+\alpha}(N)$), where $\ti \W_j$ are wedges with edge $\n$ in $\C^{k+m}$. 
 We thus apply  proposition \ref{p1} (or proposition \ref{p12}) to each
$U_j$ to obtain a holomorphic extension to wedges $\W'_j$ all in the
direction
$v$. Let $\W=\cap_{j=1}^K \W'_j$ we conclude that
the function $\sum_{j=1}^{K}
U_j$ extends holomorphically  to $\W$ and $\sum_{j=1}^{K}U_j=u$ on $\n$. 
This concludes the proof of proposition \ref{t'1}. $\blacksquare$

\bigskip

\noindent {\bf Proof of Theorem \ref{t1}.} Denote the coordinates
$(z,w',w'')$
in
$\C_z^k\times \C^m_{w'}\times \C^n_{w''}$. Recall that 
$  N=\{(  \n,h(  \n))\}$.
Consider the CR map $h: 
\n \to \C^n$. By proposition \ref{t'1}, each $h_j$ extends holomorphically
to some wedge $\W_j$ in any complex transversal direction $v$.
 Set $\W =\cap_{j=1}^n \W_j$, $\W \not = \emptyset$
since $v \in \W$.
Define $F:(  \n,0) \to (\n,\kappa h(  \n))$ where $\kappa \in \R^*$ by
$$F(z,w',w'')=(z,w',w''+\kappa h(z,w')).$$
Clearly, there exists $\kappa\not = 0$ so that on  $\overline \W$, the
Jacobian
of $F$ is non zero. Hence $F$ is a biholomorphism from $\W$
to
$F( \W)$ extending to a $\CC^{1+\alpha}$ diffeomorphism from $ \W \cup
(\n \times \{0\})$ to
$F( \W \cup
(\n \times \{0\} ))$.
Note that since $dh(0)=0$, $F$ is tangent to the identity at
the origin, hence, there exists $\W'$ a wedge in $\C^{k+m+n}$ of
direction
$v$ such that $\W' \subset F( \W)$.
 Thus any decomposable distribution (resp $\CC^{0+\alpha}$ function) on $N$ extends
holomorphically to the complex transversal wedge $\W'$.
Note then that 
the functions  $f_j=w''_j-h_j$ are holomorphic on a wedge $\W_v$
and null on $N$ and they clearly verify the desired conclusions.

\bigskip

\noindent {\bf Proof of Corollary \ref{D}.} Let $M$ and $N$ be as in
the hypothesis of the corollary. After a linear of variables, we
may assume that $p_0=0$ and that near the origin, $M$ is parametrized
by
$$M=\{z,u+iv(z,u):(z,u) \in \C^k \times \R^{p-k}\}.$$
By the implicit function theorem, we may assume that $N$ is given
as a subset of $M$ by

$$
\begin{cases}
u_{p-k-n}=\mu_1(z,u_1,...,u_{p-k-n-1}),...,
u_{p-k}=\mu_n(z,u_1,...,u_{p-k-n-1}),\\
\mu(0)=d\mu(0)=0.
\end{cases}
$$
Denote by $s=(u_1,...,u_{p-k-n-1}) \in \R^m$ and
$t=(u_{p-k-n},...,u_{p-k}) \in \R^n$. Setting $t'=t-\mu$, in the $(z,s,t')$
coordinates,
we have $N$ given as a subset of $M$ by $t'=0$ and
$$N=\{(z,w'(z,s),h(z,w'):(z,s)\in \C^k  \times \R^m\},$$
where $h$ is a CR map from $\n:=\{z,w'(z,s)\}$. We can now apply 
 theorem
\ref{t1} to obtain the CR extension as the restriction of the
holomorphic
extension of $f$ to $\W \cap M$.
 The second part of the corollary follows in the same manner.
$\blacksquare$

\section{ Non Extension Example}

We will now construct an example of an abstract CR structure $(M,{\mathcal V})$
in which there are no local CR extension property.

Set $L={{\partial} \over {\partial \bar z}} +C(z,s,t){{\partial} \over {\partial s   }}
+tD(z,s,t){{\partial} \over {\partial t   }}$ and define
$L^0=L|_{t=0}$.
Proposition \ref{P} states that
 there exist  $L$ as above and $h$ a real
analytic function, with $L^0(h)=0$ and such that the equation 
$L(h+tg)=0$ has no solution for $g \in {\mathcal C}^1$.

\bigskip

\noindent {\bf Proof of Proposition \ref{P}.} We first construct $L^0$. Let
$f: {\C} \to {\R}$ be a real analytic function such that
there exists $g \in {\mathcal C}^{\infty}_0(B_{\epsilon}(0))$ (the
neighborhood
is taken in ${\R}^3$)
where the equation 
$$L^0(u)=[{{\partial} \over {\partial \bar z}}-if_{\bar z}{{\partial} \over
{\partial s   }}](u)=g$$
 is not solvable in any neighborhood of the
origin in ${\R}^3$. (H\" ormander's Theorem)(\cite{[Ho]}p157).

\bigskip

\begin{lemma} \label{l1} There exists $\eta=\eta(z,s) \in
{\mathcal C}^{\omega}$ such that $L^0(\eta) \not = 0$ and
the equation $L^0(u)=e^{\eta}g$ is solvable nowhere.
\end{lemma}

\bigskip

\noindent {\bf Proof of the Lemma \ref{l1}.} Let $\eta \in
{\mathcal C}^{\omega}$ such that $L^0(\eta) \not = 0$. We note that if
$h \in
{\mathcal C}^{\omega}$ is such that  $L^0(h)=0$ and $h$ does not vanish in
some
neighborhood of the origin in ${\R}^3$, then one of the two
equations
is not locally solvable 
$$\begin {cases}
L^0(u)=e^{\eta}g,\\
L^0(u)=(e^{\eta}+h)g.\\
\end{cases}
$$
Indeed, if both of these equations were solvable, with solutions
$u_1$ and $u_2$, on some neighborhoods
$U_1$ and $U_2$ of the origin, then on $ U_1 \cap U_2$ we would have,
by
setting $u=u_2-u_1$ 
$$L^0(u)=L^0(h{{u}\over {h}})=hL^0({{u}\over {h}})=hg.$$
So we conclude that $L^0({{u}\over {h}})=g$, contradicting our choice of $L^0$.

\bigskip

Without loss of generality, assume $L^0(u)=(e^{\eta}+h)g$ is not locally
solvable. To finish the proof of the lemma, we wish to find $h$ such
that $h \not = 0$, $L^0(h)=0$ and $e^{\eta}+h \not = 0$, then we will
set $e^{\eta}+h=e^{\widetilde \eta}$. By Cauchy-Kovalevsky 
 we solve the equation

$$\begin {cases}
L^0(v)=0,\\
v|_{Re[z]=0}=e^{\eta}|_{Re[z]=0}.\\
\end{cases}$$

We thus have $v=e^{\eta(0,Im[z],s)}+Re[z]\zeta$ and consequently
$(e^{\eta}+v)(0,0,0) \not =0$. Therefore by eventually shrinking our
neighborhoods
we get the desired function. This completes the proof of Lemma \ref{l1}.

\bigskip

We are now ready to define $L$ from $L^0$, set 

$$L=L^0+t(g{{\partial} \over {\partial s   }}+L^0(\eta)
{{\partial} \over {\partial t   }})$$

\noindent {\bf Claim:}{\it  The function $h=s+if(z)$ admits no CR
extension to $(M,L)$}.

\bigskip

Indeed, $L^0(h)=0$. Suppose there exists $v \in {\mathcal C}^1$
such that
$L(h+tv)=0$ has a local solution, then, we note that
$L(h)=-tg$, thus we have $L(h+tv)=-tg+tL(v)+tvL^0(\eta)=0$.
So
$-g+L(v)+vL^0( \eta)=0$. Set $v=v_0(x,s)+tv_1(x,s,t)$, then
$L(v)=L^0(v_0)+tG$, thus equating terms with no $t$ and
multiplying by
$e^{\eta}$, we get
$$e^{\eta}(L^0(v_0)+v_0L^0(\eta))=L^0(v_0e^{\eta})=e^{\eta}g$$
Contradicting the Lemma \ref{l1}.

\bigskip

\noindent {\bf Remarks:} (a) There are plenty of non zero CR functions on
$(M,{\mathcal  V})$, any holomorphic function of $z$ is CR and
we can also find functions of $z$ and $t$ that are CR, indeed, 
if $f=f(z,t)$ then by
Cauchy-Kovalevsky
one can solve $L(f)=0$ with non zero Cauchy data.

\bigskip

\noindent (b) The CR structure $(M,{\mathcal  V})$ defined
above is not realizable in ${\C}^3$.

\bigskip

By realizable, we mean that there does not exist $
\Phi_1,\Phi_2,\Phi_3$, complex valued functions, such that
$L(\Phi_j)=0$ and $d\Phi_1 \wedge d\Phi_2 \wedge
d\Phi_3 \not = 0$ in a neighborhood of the origin.

\bigskip

If $(M,{\mathcal  V})$ was realizable, then any real analytic CR
function on $(N,{\mathcal  V}_0)$ would admit a CR extension to 
$M$, since $L^0$ is real analytic.

\bigskip

\noindent Nicolas Eisen\\
D\'epartement de Math\'ematiques, UMR 6086 CNRS\\
Universit\'e de Poitiers\\
eisen@math.univ-poitiers.fr

\end{document}